\documentclass[12pt]{amsart}

\usepackage{cite}

\usepackage{kotex}
\usepackage{lineno,hyperref}
\modulolinenumbers[5]
\usepackage{epstopdf}
\usepackage{pifont}
\usepackage{latexsym}
\usepackage{graphics}
\usepackage{graphicx}
\usepackage{float}
\usepackage[dvips]{xy}
\xyoption{all}
\usepackage{ifpdf}
\usepackage{multirow}
\usepackage{amssymb, amsthm, amsmath}
\usepackage{hhline}
\usepackage{centernot}
\usepackage{verbatim, comment, color}
\usepackage{enumerate}
\usepackage{setspace}
\usepackage{lipsum}
\usepackage{subcaption}

\usepackage{mathtools}

\newtheorem{theorem}{Theorem} 
\newtheorem{proposition}[theorem]{Proposition}

\newtheorem{remark}[theorem]{Remark}
\newtheorem{corollary}[theorem]{Corollary}
\newtheorem{definition}[theorem]{Definition}
\newtheorem{example}[theorem]{Example}
\newtheorem{claim}[theorem]{Claim}

\newtheorem{question}[theorem]{Question}

\newcommand{\ba}{\begin{align}}
\newcommand{\ea}{\end{align}}  

\newcommand{\be}{\begin{equation}}

\newcommand{\ee}{\end{equation}}
\newcommand{\bea}{\begin{eqnarray}}
\newcommand{\eea}{\end{eqnarray}}
\newcommand{\barr}{\begin{array}}
\newcommand{\earr}{\end{array}}
\newcommand{\bn}{\begin{enumerate}}
\newcommand{\en}{\end{enumerate}}
\newcommand{\bi}{\begin{itemize}}
\newcommand{\ei}{\end{itemize}}
\newcommand{\bbbm}{\begin{pmatrix}}
\newcommand{\eeem}{\end{pmatrix}}

\newcommand{\bbN}{{\bf N}}

\newcommand{\bbS}{{\bf S}}

\newcommand{\co}{{c_d}}
\newcommand{\cf}{{C_d}}

\newcommand{\R}{{\mathbf R}}

\newcommand{\al}{\alpha}

\newcommand{\la}{\lambda}

\newcommand{\tta}{\theta}
\newcommand{\ignore}[1]{}{}
\newcommand{\noin}{\noindent}

\newcommand{\nn}{\nonumber}

\newcommand{\q}{\quad}

\newcommand{{\QED}}{{\hfill QED} \smallskip}

\renewcommand{\subset}{\subseteq}

\renewcommand{\phi}{\varphi}
\newcommand{\cal}{\mathcal}

\DeclareMathOperator{\conv}{conv}
  \DeclareMathOperator*{\diam}{diam}

\numberwithin{equation}{section}
\numberwithin{theorem}{section}
    
\begin{document}
\title[Angle bounds limiting the cardinality of sets in ${\bf R}^d$]
{On the cardinality of sets in ${\bf R}^d$ obeying a 
slightly obtuse angle bound$^*$} 
\thanks{$^*$\em{The authors thank Sungjin Kim and Bal\'asz Gerencs\'er for fruitful remarks, and
Dmitry Zakharov for directing their attention to his work with Kupavskii \cite{KupavskiiZakharov20} and its connection to \cite{EF83}.}
TL is grateful for the support of ShanghaiTech University, and in addition, to the University of Toronto and its Fields Institute for the Mathematical
Sciences, where parts of this work were performed.  RM  acknowledges partial support of his research by  the Canada Research Chairs program and
Natural Sciences and Engineering Research Council of Canada Grant 2015-04383 and 2020-04162.
\copyright 2020 by the authors.}

\date{\today}

\author{Tongseok Lim and Robert J. McCann}
\address{Tongseok Lim: Krannert School of Management \newline  Purdue University, West Lafayette, Indiana 47907}
\email{lim336@purdue.edu / tlim0213@outlook.com}
\address{Robert J. McCann: Department of Mathematics \newline University of Toronto, Toronto ON Canada}
\email{mccann@math.toronto.edu}

\begin{abstract}  In this paper we explicitly estimate the number of points in a subset  $A \subset \R^{d}$
as a function of the maximum angle $\angle A$ that any three of these points form,  
provided $\angle A < \theta_d := \arccos(-\frac 1 {d}) \in (\pi/2,\pi)$. We also show $\angle A < \theta_d$ ensures
that $A$ coincides with the vertex set of a convex polytope.
This study is motivated by a question of Paul Erd\H{o}s and indirectly by a conjecture of L\'aszl\'o Fejes T\'oth.
\end{abstract}

\maketitle

\section{Introduction}

Let us begin with a simple definition.
\begin{definition}[Angle bound] For $A \subset \R^{d}$ let $\angle A$ denote the smallest value of $\tta \in [0, \pi]$, 
such that no triple of points $x,y,z \in A$ determine an angle  $\angle xyz$ greater than $\theta$, i.e. 
$\frac{x-y}{|x-y|} \cdot \frac{z-y}{|z-y|} \ge \cos \theta$ for all $x,y,z \in A$.
\end{definition}

For a set $A$, we denote by $|A| \in \{0,1,..., \infty\}$ its cardinality.  
In the 1950s, Paul Erd\H{o}s raised the following conjecture 
\cite{Erdos57}:

\medskip

{\bf Conjecture:} If $A \subset \R^{d}$ satisfies $\angle A \le \pi/2$, then $|A| \le 2^{d}$. 
\medskip

The conjectured bound $2^{d}$ is obviously sharp, being achieved by the vertices of the 
hypercube in $\R^d$. 
For $d=3$ he had advertised the problem a decade earlier \cite{Erdos48},  which he claimed \cite{Erdos57} lead to an unpublished
solution by Kuiper and Boerdijk.
For $d \ge 3$ the conjecture was resolved affirmatively by Danzer and Gr{\"u}nbaum \cite{DG62}, who 
established it through a chain of remarkable inequalities reproduced in \cite{AZ18}.
They also asked the following natural question: if the angle bound $\tta$ is acute, i.e. $\tta \in [0, \frac\pi 2]$, what would be the optimal upper bound for $|A|$ subject to 
 the strict inequality $\angle A < \tta$? Danzer and Gr{\"u}nbaum had raised a conjecture on this question which was disproved two decades later by Erd\H{o}s and F{\"u}redi \cite{EF83}. Since then the question has remained both interesting and challenging, as e.g. Gerencs{\'e}r and Harangi \cite{GH19}, \cite{GH19-2} and 
Aigner and Ziegler \cite{AZ18} discuss; a significant stream of research has focused on how large a set $A$ can be while satisfying the strictly acute bound  $\angle A < \frac\pi 2$.\\

Now what if the angle bound $\tta$ is obtuse, i.e. $\tta \in (\pi/2, \pi)$? {In this case much less appears to be known. It is known however that the set $A \subset \R^d$ subject to $\angle A \le \tta$ must be finite; see e.g. Otsetarova \cite{O17}. 
In the plane, the number of obtuse angles determined by a set is known to grow as a cubic function of its cardinality \cite{CCEG79} \cite{F-MHT14}.
Here we address the following natural question: 


\begin{question}\label{Q:N} 
Given $\tta \in (\pi/2, \pi)$, estimate the smallest  $N=N(\tta, d) \in \bbN$ 
(if any exists) such that for all $A \subset \R^{d}$ satisfying $\angle A \le \tta$, we have $|A| \le N$?
\end{question}


Since $N(\pi,d)=+\infty$,  it is not obvious that $\tta<\pi$ makes $N(\tta,d)$ finite;
we describe in an appendix below how this fact can be inferred by refinining the conclusion of Erd\H{o}s and F\"uredi \cite[Theorem 4.3]{EF83}, 
who used sphere packing and covering asymptotics to derive 
two-sided bounds on $N(\tta,d)$ for $\tta$ near $\pi$;  
their statement is predicated on the unstated requirement that our $N$ (their $n$) be sufficiently large; we clarify how large in Theorem \ref{T:E-F}.
After the current manuscript was first posted to the arXiv, 
we also learned that Kupavskii and Zakharov \cite{KupavskiiZakharov20} had used  a method inspired by Erd\"os and F\"uredi's to show for each 
$\frac\pi2 <\tta < \pi$ the growth of $N(\theta,d)$ to be {\em doubly} exponential in $d$, but with $\theta$ dependent rates which are not estimated explicitly.
In this paper, we use an entirely different approach to give the first explicit upper bound for $N(\theta,d)$ in the region $\theta<\theta_d$, when
\begin{equation}\label{thetad}
\theta_d := \arccos\Big(-\frac{1}{d}\Big) \q \in \q (\frac\pi 2,\pi)
\end{equation}
shrinks to $\pi/2$ as $d \to \infty$.
 Instead of relying on the sphere packing and covering arguments (see Appendix A for more details),  we observe that $\theta:=\angle A < \theta_d$ 
implies the points of $A$ form the vertices of a convex polytope.  Our bound for the number of vertices of this polytope in terms of $\theta$ and $d$ then 
follows from the Gauss-Bonnet theorem. This argument is a higher dimensional version of the planar statement that: if all angles of a convex polygon are at most $2\pi/3$, then the polygon has at
most six vertices.

In the first circulated draft of our companion work \cite{LM20_2}, we exploited the present result to attack a conjecture of Fejes T\'oth,
concerning the placement of a large number of lines through the origin so as to maximize the
expected acute angle between them.  There instead we were interested in maximizing a power $\alpha>1$ of the (renormalized) angle between each pair of lines.
Having already established the analogous conjecture in the limiting case $\al=\infty$ \cite{LM21},  we used the present results to extend this conclusion to large finite values of $\al$.
In subsequent drafts of \cite{LM20_2} however, this argument has been replaced by a different approach (based in part on an appendix authored by Bilyk, Glazyrin, Matzke, Park and Vlasiuk) which allows us to extend our conclusion to the larger range $\al>\al_{\Delta^d}$, where $\al_{\Delta^d}<2$ in certain cases. 

Changing the ambient space from $\R^d$ to $\R^{d+1}$ hereafter, let $\cal H^{d}$ be the area measure on $\bbS^{d}=\{x \in \R^{d+1} \ | \ |x|=1\}$, 
and let  $\omega_{d}=\cal H^{d}(\bbS^{d})$ denote the total area of $\bbS^{d}$, e.g. $\omega_1=2\pi$, $ \omega_2=4\pi$. For $\eta \in (0,\frac{\pi}{2})$,  let $f_{d}(\eta)$ 
denote the fraction of $\bbS^{d}$ covered by generalized normals to the cone
\begin{equation}\label{ueta cone}
C_{u, \eta} :=\{x \in \R^{d+1} \ | \ u \cdot x \ge |u| |x|\cos \eta \}
\end{equation}
of half-angle $\eta$ around $0\ne u \in\R^{d+1}$, so that
\begin{align}
f_{d}(\eta) 
&:= \frac1{\omega_{d}} \cal H^{d} (  \{x \in \bbS^{d} \ | \ x \cdot z \le 0 \, \mbox{ for every }  z \in C_{u,\eta}\} )
\nonumber \\ &= \frac{\omega_{d-1}}{\omega_{d}} \int_0^{\frac\pi 2-\eta} \sin^{d-1}(t) dt
\label{cone curvature} \\ 
\nonumber &= \Big( \int_0^{\pi} \sin^{d-1}(t) dt\Big)^{-1} \int_0^{\frac\pi 2-\eta} \sin^{d-1}(t) dt.
\end{align}
Equivalently, $f_{d}(\eta)$ is the zeroth curvature measure that the cone $C_{u,\eta}$ assigns to its vertex,
in the terminology of Federer \cite{Federer59} \cite{Schneider78}.

Our result  is the following.
\begin{theorem}[Cardinality under possibly obtuse angle bound]
\label{boundby91}
Fix $d \ge 1$ and $ 0<\tta <\theta_{d+1}$ and
$\eta_{d}(\theta):= \arcsin \Big(\frac{\sin(\theta/2)}{\sin(\theta_{d}/2)} \Big) \in (0,\frac\pi 2)$.
If $A \subset \R^{d+1}$ satisfies $\angle A \le \theta$, then \eqref{thetad} and \eqref{cone curvature} yield
\begin{equation}\label{cardinality bound}
|A|      \le   \frac{1}{f_{d}(\eta_{d}(\theta))}.
\end{equation}
\end{theorem}

\begin{remark}[Dimensional monotonicity]\label{R:estimates}
For fixed $\theta$,  our cardinality bound $1/f_{d}(\eta_{d}(\theta))$ increases with dimension since
 $f_{d}(\eta)$ is monotone decreasing with respect to both $\eta$ and $d$:  indeed the difference of 
averages
 \begin{align*}
 \frac{\partial \log f_{d}(\eta)}{\partial d} 
 =&  
 \Big( \int_0^{\frac \pi 2 -\eta} \sin^{d-1}(t) dt\Big)^{-1} \int_0^{\frac\pi 2-\eta} ( \log \sin t) \sin^{d-1}(t) dt
\\ &-  \Big( \int_0^{\pi} \sin^{d-1}(t) dt\Big)^{-1} \int_0^{\pi} (\log \sin t)\sin^{d-1}(t) dt
 \end{align*}
 is negative since $t \in [0, \pi] \mapsto \log\sin t $ is symmetric about $\pi/2$ and increasing on $[0, \pi/2]$.
\end{remark}

\begin{example}[Explicit bounds in low dimensions]
If $d=1$ then $\eta_1(\theta)=\frac\theta{2}$ and $f_1(\eta)=\frac12-\frac\eta\pi$ so $f_1(\eta_1(\theta))^{-1}= \frac{2\pi}{\pi-\theta}$ yields the known sharp values $\{2,3,4,\infty\}$ corresponding to $\theta \in \{0,\frac\pi 3, \frac \pi 2, \pi\}$.
However $f_2(\eta) =\frac12(1-\cos(\frac\pi 2 - \eta))$
yields a bound $f_2(\eta_2(\frac\pi 2))^{-1} \approx 10.9$
worse than the sharp value $8$ attained by the vertices of the cube in $\R^3$.
\end{example}

For non-obtuse sets, the bound \eqref{cardinality bound} becomes less and less accurate with increasing dimension,
growing slightly faster than exponentially:

\begin{claim}[Asymptotic bounds as dimension increases]
\label{claim}
\[
\ f_d(\eta_d(\frac\pi 2))^{-1}
<2 (\frac \pi 2)^{2d-1}e^{(d\log d)/2}[1 + O(d^{-1})].
\]
\end{claim}

This growth, while faster than the known sharp value of $2^{d+1}$, is slow 
compared to the doubly exponential growth of $N(\theta,d)$ observed for each fixed 
$\tta \in (\frac\pi2,\pi)$ in
\cite{KupavskiiZakharov20}, suggesting a rich range of intermediate asymptotic behaviour
in the narrowing region $\frac\pi 2< \theta < \theta_d$ of the $(\theta,d)$ plane.
It is also smaller that the doubly exponential upper bound obtained in Appendix \ref{S:Refined}
using Erd\H{o}s-F\"uredi's sphere-covering arguments (with $\cf$ independent of $d$),  which confirms that our bounds represent a significant improvement on the 
state of the the art at least in their limited range of validity.

\bigskip
\noin{{\bf Proof of Claim \ref{claim}.} From the identity
\[
\sin(\theta_d/2) :=   \sqrt{\frac{d+1}{2d}}
\]
which follows from \eqref{thetad} (and which is reasserted below in Theorem~\ref{T:Dekster}) we derive

\begin{align*}
\eta_{d}(\frac\pi2)
&=  \arcsin \Big(\sqrt{1 - \frac1{d+1}}\Big)
\\ & <  \frac\pi 2 -\sqrt{2(1-\sqrt{1-(d+1)^{-1}})}
\\& = \frac\pi 2 -d^{-1/2}(1+O(d^{-1}))
\end{align*}
where $\arcsin(t) < \frac\pi 2 -\sqrt{2-2t}$ has been used.  Estimating \eqref{cone curvature}
using $\frac2\pi  < t^{-1}\sin(t) < 1$ 
on $t \in (0,\frac\pi 2)$ yields
\begin{align*}
f_{d}(\eta) 
& >  \Big( 2\int_0^{\pi/2} t^{d-1} dt\Big)^{-1} \int_0^{\frac\pi 2-\eta} (\frac{2}{\pi} t)^{d-1} dt
\end{align*}
whence
\begin{align*}
 f_d(\eta_d(\frac\pi 2))^{-1}
 &< 2 (\frac \pi 2)^{2d-1}d^{d/2}[1 + O(d^{-1})]
\end{align*}
as desired.
\QED

In the next section we will prove Theorem \ref{boundby91} using two propositions
which may have independent interest.  
Proposition \ref{simplexbound} shows that the strict inequality $\angle A <\theta_{d}$ cannot hold unless no point in 
$A\subset \R^{d}$ is a convex combination of $d+1$ others.  On the other hand,  if no point in $A$ is a 
 convex combination of other points,  then $A$ consists precisely of the vertices of a convex polytope.  
 Proposition \ref{polytopebound} combines a spherical diameter-to-radius inequality \cite{D95} with the generalized Gauss-Bonnet
 theorem to estimate $|A|$ in terms of the angle bound $\angle A \le \theta$ in this case.
 {In an appendix we review the sphere-packing and covering arguments of Erd\H{o}s and F\"uredi for comparison,  and show how they can be extended
 to sets which need not be too large.}

\section{Proofs} 
 Let $\conv(A)$ denote the convex hull and $\mathop{\rm int}(A)$ the interior of any 
 subset $A\subset \R^d$.
 
\begin{proposition}[Angle estimates from an interior point of a simplex]
\label{simplexbound}
Let $d \ge 2$ and let $\{v_0,...,v_d\} \subset \R^d \setminus \{0\}$ be vertices of a $d$-dimensional simplex containing the origin. Let $w_i=  \frac{v_i}{|v_i|}$. Then 
\be \label{dotbound}
\min_{0 \le i < j \le d} w_i \cdot w_j \le -\frac{1}{d}
\ee
and equality holds if and only if $\conv\{w_0,...,w_d\}$ is a regular simplex. 
\end{proposition}
\noin{\bf Proof.} The proposition clearly holds for $d=2$: 
since at least one of the three angles at the origin
must exceed $2\pi/3$ unless $\{w_0,w_1,w_2\}$ is equilateral. We will proceed by an induction on dimension $d$. 
Let $V= \{v_0,...,v_d\}$, so that $\conv (V)$ is the $d$-dimensional simplex. If the origin lies on the boundary of $\conv (V)$, then the induction hypothesis yields $\min_{i \ne j} w_i \cdot w_j \le -\frac{1}{d-1}$. So let us assume $0 \in {\rm int} (\conv (V))$. Let $w_i = \frac{v_i}{|v_i|}$. We claim $W= \{w_0,...,w_d\}$ also forms vertices of a $d$-dimensional simplex containing the origin in its interior. To see this, observe 
\begin{align*}
&0 \in {\rm int}( \conv (V)) \\
&\iff \text{for every } i, \ v_i = - \sum_{j \ne i} \la_j v_j \text{ for some } \la_j >0  \\
& \iff  \text{for every } i, \ w_i = - \sum_{j \ne i} \la_j w_j \text{ for some } \la_j >0 \\
&  \iff 0 \in {\rm int} (\conv (W))
\end{align*}
which proves the claim, and gives $a_i >0$ such that $\sum_{i=0}^d a_i w_i =0$. 
Without loss of generality assume $a_0 = \min_i \{a_i\}$, $a_1= \max_i\{a_i\}$. Set $b_i=\frac{a_i}{a_0}$, so that $w_0 = -\sum_{i=1}^d b_iw_i$. Now we claim:
\be
\min_{1 \le i < j \le d} w_i \cdot w_j \ge -\frac{1}{d} \ \text{ implies } \ w_0 \cdot w_1 \le -\frac{1}{d}. \nn
\ee
To see this, observe
\begin{align*}
w_0 \cdot w_1 &=  -\sum_{i=1}^d b_i w_1 \cdot w_i  \\
& \le -b_1 + \frac{1}{d}  \sum_{i=2}^d b_i  \\
&\le -b_1+ \frac{1}{d}  \sum_{i=2}^d b_1 =   -\frac{b_1}{d} \\
& \le -\frac{1}{d}
\end{align*}
since $b_1 \ge 1$. We have shown \eqref{dotbound}, and in view of the above inequalities, we see that equality holds in \eqref{dotbound} if and only if $w_i \cdot w_j =-\frac{1}{d}$ for all $1 \le i < j \le d$ and $b_1=1$, that is $\sum_{i=0}^d w_i=0$. By taking dot product with $w_j$, $j \ne 0$ on the last identity, we conclude that equality holds in \eqref{dotbound} if and only if $w_i \cdot w_j =-\frac{1}{d}$ for all $0 \le i < j \le d$, that is, $\conv\{w_0,...,w_d\}$ forms a regular $d$-dimensional simplex. 
\QED

\begin{corollary}[Deciding when points lie in convex position]
\label{coro}
Let $d \ge 2$ and $A \subset \R^d$. If $\angle A < \theta_d := \arccos(-\frac1 d)$, then $A$ consists of the vertices of a convex polytope (not necessarily $d$-dimensional).
\end{corollary}

\noin{\bf Proof.} If not, there exist $v \in A$ and $\{v_0,...,v_k\} \subset A \setminus \{v\}$ which forms vertices of a $k$-dimensional simplex such that $v \in \conv\{v_0,...,v_k\}$ by Carath{\`e}odory's theorem. Proposition \ref{simplexbound} 
then yields $i,j$ such that $\angle (v_i,v,v_j) \ge \theta_k \ge\theta_d$, a contradiction. \QED
\\

To derive the desired cardinality bound for convex polytopes,  we use the following spherical
 version of Jung's Theorem \cite{ Jung01}
relating diameter to radius bounds in flat space, established on the unit sphere $\bbS^{d}$ equipped with its standard round metric
by Dekster   \cite[Theorem 2]{D95}:

\begin{theorem}[Dekster, 1995]\label{T:Dekster}
Let $R \in (0, \frac{\pi}{2})$. If $H \subset \bbS^{d}$ satisfies $\diam H \le 2\arcsin\big(\sin(\theta_{d}/2) \sin R \big) \in (0,\pi)$, then $H$ can be contained in a closed ball of radius $R$ in $\bbS^{d}$.
Here $\sin(\theta_{d}/2) = \sqrt{\frac{d+1}{2d}}$ from \eqref{thetad}.
\end{theorem}

\begin{proposition}[Cardinality bound for convex polytopes]
\label{polytopebound} Let $d \ge 1$ and $0<\theta <  \theta_{d}$. If  $V=\{v_1,...,v_n\}$ consists of the vertices of a convex polytope $\conv(V) \subset \R^{d+1}$ with non-empty interior and satisfies  $\angle V \le \theta$,  then $n \le  \frac{1}{f_{d}(\eta)} $ with $\eta= \eta_{d}(\theta)=\arcsin\Big( \frac{\sin (\theta/2)}{\sin(\theta_{d}/2)}\Big) \in (0,\frac\pi 2)$.
\end{proposition}

\noin{\bf Proof.}
\noin For $\kappa \in (0, \frac{\pi}{2})$ and $u,v \in \R^{d+1}$, 
use \eqref{ueta cone} to 
define a family of cones $v+C_{u,\kappa} \subset \R^{d+1}$ with vertex $v$, direction $u \ne 0$, and half-angle $\kappa$. 
 Now consider any of the given vertices $v_i$; by translating the polytope we may assume $v_i=0$. By dilating $\la V$ with sufficiently large $\la >0$ we also assume that only those edges emanating from $v_i=0$ meet $\bbS^{d}$, 
 and set $H \subset \bbS^{d}$ to be the intersection points of those edges with $\bbS^{d}$. Notice $\angle V \le \tta$ implies $ \diam H \le \tta$, hence Dekster's theorem implies, with $\eta$ in place of $R$, that $H$ (and hence $V$) is contained in the cone $v_i+C_{u_i,\eta}$ for some unit vector $u_i$. Since the conclusion is invariant under translations and dilations, we deduce there exists unit vectors $u_1,...,u_n$ such that
\be
\label{covering}
\conv V \subset \bigcap_{i=1}^n v_i+C_{u_i, \eta}.
\ee
Let $f_i$ denote the fraction of $\bbS^{d}$ occupied by the generalized normals to $\conv(V)$ at $v_i$:
$$
f_i := \frac1{\omega_{d}}\cal H^{d} ( \{x \in \bbS^{d} \ | \ (x-v_i) \cdot (z-v_i) \le 0 \, \mbox{ for every }  z \in \conv V \} ).
$$
Then 
$$
\sum f_i =1,
$$
an intuitive fact which can also be seen as a consequence of, e.g., Federer's generalization of the Gauss-Bonnet formula, 
which asserts that for any convex body the zeroth curvature measure of the entire convex body is unity,
coinciding with its Euler-Poincar\'e characteristic  \cite{Federer59};  
in the case of a convex polytope $\conv(V)$, 
the zeroth curvature measure vanishes except at the vertices of the body,  and assigns mass $f_i$ to $v_i$.

Now by the covering property \eqref{covering}, we have
\be
\label{curvature1}
1=\sum_{i=1}^n f_i \ge \sum_{i=1}^n f_{d}(\eta) = n f_{d}(\eta),
\ee
whence $n \le 1 / f_{d}(\eta)$.
\QED
\\

\noin{\bf Proof of Theorem \ref{boundby91}.} Let $A \subset \R^{d+1}$ satisfy $\angle A \le \tta <  \theta_{d+1}$. We may assume that $A$ is not contained in any $d$-dimensional hyperplane, since otherwise we may apply an induction on dimension
using the monotonicity of $d \mapsto f_{d}(\eta_{d}(\theta))$ established in Remark \ref{R:estimates}. Then by Corollary \ref{coro}, $A$ is the set of vertices of a $(d+1)$-dimensional convex polytope, and the theorem follows from Proposition \ref{polytopebound} and $\tta_{d+1}<\tta_d$.
\QED

\appendix

{
\section{Finite-size Erd\H{o}s and F\"uredi bounds} 
\label{S:Refined}

 Erd\H{o}s and F\"uredi  relate the quantity we have estimated 
to sphere-packing and sphere-covering bounds on the unit sphere $\bbS^{d-1}$ 
in \cite[Theorem 4.3]{EF83}.
Since their published statement
contains at least one misprint  (and some tacit hypotheses, as observed in \cite{KupavskiiZakharov20}), let us recount 
their argument and refine their conclusions, so that they apply to all sets, and not only to large ones.
Using our angle bound
$$
\angle P := \sup_{x,y,z \in P} \angle xyz,
$$ 
let
$$
\alpha_d(n) := \inf_{P \subset \R^d \atop |P|=n} \angle P.
$$
denote the minimal maximal angle made by $n$ points in $\R^d$,
e.g. $\al_d(3)=\pi/3=60^\circ$,  $\al_3(4)\approx 109.5^\circ$, etc.
Note $P\subset P'$ implies $\angle P \le \angle P'$, yielding
monotonicity of this minimax $\al_d(n) \le \al_d(n+1)$ with respect to the size of the sets being considered.   This means 
$$
N_d(\tta) := \sup_{\al_d(n) \le \tta} n
$$
is a non-decreasing inverse to $\al_d(n)$ (and in fact agrees with $N(\theta,d)$ from our Question \ref{Q:N}).
Any upper (respectively lower)
bound for $\alpha_d(n)$ can therefore be translated into a lower (respectively upper) bound for 
$N_d(\tta)$, and vice versa.  
Theorem 4.3 of \cite{EF83} asserts
\be\label{EFerror}
\frac\pi{(\log_2 n)^{1/(d-1)}}
< \pi - \al_d(n) <
\frac{4\pi}{(\log_2 n)^{1/(d-1)}}
 \ee
under the tacit assumption that 
$n$ needs to be sufficiently large.
 Inspection of the proof suggests the intended statement should perhaps have been
 that 
\begin{equation}\label{revised}
 \frac{\co}{(1+\log_2 n)^{1/(d-1)}}
<  \pi - \al_d(n) 
< \frac{\cf}{(-1+\log_2 |n-1|)^{1/(d-1)}},
\end{equation}
where $\co$ and $\cf$ are certain constants independent of $d$, to be estimated presently,
and the lower bound is only asserted for large $n$.  
Note \eqref{revised} is equivalent to the bound
\be\label{KZ}
2^{(\co/(\pi-\tta))^{d-1}-1} < N_d(\theta) < 1+  2^{(\cf/(\pi-\tta))^{d-1}+1}
\ee
while the requirement that $N_d(\theta)$ be large translates to $\tta$ being close to $\pi$.  
For $\theta>\pi-\co$ the bounds \eqref{KZ} grow doubly exponentially with dimension $d$.
Kupavskii and Zakharov established doubly exponential bounds for $N_d(\theta)$ in the full range $\frac\pi2< \tta < \pi$,
but their rate constants are not explicit and must evidently depend on $\tta$ near $\frac\pi2$ \cite{KupavskiiZakharov20}.

Let us now recall the arguments which establish \eqref{revised}.
They begin with two crucial statements concerning the packing and covering of the unit sphere in $\R^d$
with balls of radius $\rho/2$: namely, that there exist 
positive functions $\co(\rho)$ and $\cf(\rho)$ such that:

{\bf (Packing)}: For all $0<\rho < \bar\rho< \frac\pi 2$, 
there exist more than $(\co(\bar\rho)/\rho)^{d-1}$ (undirected) lines through the origin  in $\R^d$  such that
any two of them determine an angle greater than $\rho$.

{\bf (Covering)}: For all  $0<\rho< \bar \rho < \pi$, 
there exist fewer than $(\cf(\bar \rho)/\rho)^{d-1}$ (undirected) lines through the origin  in $\R^d$ such that any line determines an angle less than $\rho/2$ 
with at least one of them.

Clearly it costs no generality to suppose $\cf(\rho)$ and $-\co(\rho)$ to be non-decreasing,  and to choose them to have finite non-zero limits
$\cf(0^+)$ and $-\co(0^+)$.  We may take them to be independent of dimension, or we may take them to be minimal, i.e.~to be the actual dimension and radius dependent 
packing and covering profiles of $\bbS^{d-1}$.
Erd\H{o}s and F\"uredi assert more, namely that we can take $\co(\rho)=1$ and $\cf(\rho)=4$. 
However since they give neither reference or proof we have not attempted to confirm these 
precise values, but note that they can only hold
under a smallness hypothesis on $\bar \rho$ which they presumably had in mind, since {\bf (Covering)} contradicts $\cf(\rho) =4$
when $\rho \lesssim \frac\pi 2$ if $d=2$.  
In any case it makes no difference to the remainder of
their argument.  Indeed, if we accept these two facts,  the proof of the lower bound is a simple iterative construction
while the upper bound relies on an elementary lemma from graph theory. We now recall this to refine the conclusions of 
\cite{EF83}.  See \cite{EF83} or \cite{KupavskiiZakharov20} for a fuller development of other parts of the proof.

\begin{theorem}[Finite size Erd\H{o}s-F\"uredi type bounds]\label{T:E-F}
Fix $d \ge 2$, $n \ge 3$ and $0< \bar\rho< \pi$.  If $\cf:[0,\infty)\to(0,\infty)$ 
satisfies {\rm \bf (Covering)},
then the second inequality in \eqref{revised} holds with $\cf=\cf(\bar \rho)$.
If $\bar\rho< \frac\pi 2$ and $\co:[0,\infty)\to(0,\infty)$ 
satisfies {\rm \bf (Packing)},
then the first inequality in \eqref{revised} holds with $\co=\co(\bar \rho)$ provided $n > (d-1)\log_2 \frac{\co(\bar \rho)}{\bar \rho}$.
\end{theorem}

\bigskip\noin{\bf Proof:}  
{\bf Packing bound}:  Given $0<\rho < \bar\rho< \frac\pi 2$, 
fix $m>(\co(\bar \rho)/\rho)^{d-1}$ 
and lines 
$L_1,\ldots, L_m$ such that any two of them make an angle larger than $\rho$.   If we regard them instead
as directed lines,  the angle between any two of them is less than $\pi-\rho$.  Take two points $A,B$
parallel to $L_1$, and then translate them parallel to $L_2$ far enough that the any vector from the 
original to the translated pair almost parallels $L_2$.  The largest angle between these four points,
is at worst the directed angle between $L_1$ and $L_2$, hence less than $\pi-\rho$.  Now iterate
the construction to obtain $2^m$ points having largest angle less than $\pi-\rho$.  This shows
\[
\al_d(2^m) \le \pi - \rho< \pi - \frac{\co(\bar \rho)} {m^{1/(d-1)}},
\]
and the same bound extends to $\al_d(n)$ for all $n\le 2^m$. Thus we have shown the first inequality in \eqref{revised} holds with $\co=\co(\bar \rho)$
as soon as $n > (d-1)\log_2 \frac{\co(\bar \rho)}{\bar \rho}$.

{\bf Covering bound}:  Given $0<\rho < \bar\rho< \pi$, fix $m<(\cf(\bar\rho)/\rho)^{d-1}$ 
and lines 
$L_1,\ldots, L_m$ such that any other line makes an angle less than $\rho/2$ with one of them.
For any set $P$ with $2^m+1$ points,  color the complete graph on $P$ by assigning to each
edge $xy$ a color $i$ such that $xy$ is within angle $\rho/2$ of $L_i$.  Then a lemma of Szekeres \cite{Szekeres41}
asserts the existence of an odd cycle in some color.  
(Alternatively, the last paragraph of the proof of Kupavskii and Zakharov 's Proposition 1 \cite{KupavskiiZakharov20} 
 gives a self-contained proof of an even simpler statement which also suffices for the present purpose.)
This cycle contains points $x,y,z$ 
such that all edges of the triangle $\Delta xyz$ are within angle $\rho/2$ of the same $L_i$.  The largest
angle in this triangle is at least $\pi-\rho$.  Thus $n \ge 2^m+1$ implies
$$
\al_d(n) \ge \al_d(2^m+1) \ge \pi - \rho>\pi -\frac{\cf(\bar \rho)}{m^{1/(d-1)}}.
$$
Now letting $\rho \to 0$ yields the second inequality in \eqref{revised} holds with $\cf=\cf(\bar \rho)$  for all $m$ and hence all $n$.
\QED
}

\end{document}